\documentstyle[12pt]{article}
   \font\tenmsb=msbm10 scaled\magstep 1
   \font\sevenmsb=msbm7 scaled \magstep 1
   \font\fivemsb=msbm5 scaled \magstep 1
\newfam\msbfam
      \textfont\msbfam=\tenmsb
      \scriptfont\msbfam=\sevenmsb
      \scriptscriptfont\msbfam=\fivemsb
\def\Bbb#1{{\fam\msbfam #1}}

\font\tenblah=msbm10 scaled\magstep 3
   \font\sevenblah=msbm7 scaled \magstep 1
   \font\fiveblah=msbm5 scaled \magstep 1
\newfam\blahfam
      \textfont\blahfam=\tenblah
      \scriptfont\blahfam=\sevenblah
      \scriptscriptfont\blahfam=\fiveblah

\font\tengothic=eufm10 scaled\magstep 1
\font\sevengothic=eufm7 scaled\magstep 1
\newfam\gothicfam
      \textfont\gothicfam=\tengothic
      \scriptfont\gothicfam=\sevengothic
\def\goth#1{{\fam\gothicfam #1}}

\newcommand{\mapdown}[1] {\downarrow
           \rlap{$\vcenter{\hbox{$\scriptstyle#1$}}$}  }
\newcommand{\m}{{\goth m}}

\newcommand{\verylong}{\hbox{$\hbox to .45in{\rightarrowfill}$}  }
\newcommand{\kindalong}{\hbox{$\hbox to .35in{\rightarrowfill}$}  }

\newtheorem{thm}{Theorem}[section]
\newtheorem{prop}[thm]{Proposition}

\newtheorem{cor}[thm]{Corollary}

\newtheorem{defn}[thm]{Definition}
\newtheorem{rmk}[thm]{Remark}
\newtheorem{exemp}[thm]{Example}

\newcommand{\qed}{\hskip 1cm $\rlap{$\sqcap$}\sqcup$}

\newenvironment{proof}{{\em Proof: }}{\qed\bigskip}
\newenvironment{example}{\begin{exemp}\em}{\end{exemp}}
\newenvironment{remark}{\begin{rmk}\em}{\end{rmk}}
\newenvironment{definition}{\begin{defn}\em}{\end{defn}}

\newcommand{\rank}{\mathop{\rm rank\,}}
\newcommand{\depth}{\mathop{\rm depth\,}}

\newcommand{\Ext}{\mathop{\rm Ext}\nolimits}
\newcommand{\Tor}{\mathop{\rm Tor}\nolimits}

\newcommand{\ann}{\mathop{\rm ann\,}}   
\newcommand{\im}{\mathop{\rm im\,}}
\newcommand{\coker}{{\mathop{\rm coker\,}}}

\newcommand{\Pn}{{\Bbb P^n}}

\newcommand{\projdim}{\mathop{\rm pd\,}}

\title{Submodules of the Deficiency Modules and an \\
	Extension of Dubreil's Theorem
	\footnotetext{{\em {\rm 1991} Mathematics Subject Classification.}
		14M06, 13C40}}

\author{Heath M. Martin \and Juan C. Migliore}
\date{}

\begin{document}
\maketitle
\setlength{\baselineskip}{18pt}

In this paper, we consider a basic question in commutative algebra:
if $I$ and $J$ are ideals of a commutative ring $S$, when
does $IJ = I \cap J$?  More precisely, our setting will be
in a polynomial ring $k[x_0, \dots, x_n]$, and the ideals
$I$ and $J$ define subschemes of the projective
space $\Bbb P_k^n$ over $k$.  In this situation, we are able to
relate the equality of product and intersection to the behavior
of the cohomology modules of the subschemes defined by $I$ and $J$.
By doing this, we are able to prove several general algebraic
results about the defining ideals of certain subschemes of
projective space.

Our main technique in this paper is a study of the
{\em deficiency modules} of a subscheme $V$ of $\Pn$.  These
modules are important algebraic invariants of $V$, and
reflect many of the properties of $V$, both geometric and
algebraic.  For instance, when $V$
equidimensional and $\dim V \ge 1$, the deficiency modules of $V$ are
invariant (up to a shift in grading) along the even liaison
class of $V$ (\cite{Rao}, \cite{Mig:liaison},  \cite{Sch:notes}, \cite{Hart}),
although they do not in general completely determine the even
liaison class, except in the case of curves in $\Bbb P^3$, \cite{Rao}.
On the algebraic side, at least for curves in $\Bbb P^3$, the
deficiency modules have been shown to have connections to the
number and degrees of generators of the saturated ideal
defining $V$, \cite{Mig:submodules}.  One of our main goals
in this paper is to extend
these results to subschemes of arbitrary codimension in any
projective space $\Bbb P^n$.

We now describe the contents of this paper more precisely.  In
the first section, we set up our notation and give the basic
definitions which we will use.  Then we prove our main technical
result:  if $I$ and $J$ define subschemes $V$ and $Y$, respectively,
of ${\Bbb P}^n$, we relate the quotient module $(I \cap J) / IJ$ to
the cohomology of $V$, at least when
$V$ and $Y$ meet properly.  We are then
able to give a different proof of a general statement
due to Serre about when there is an equality
of intersection and product.

In the second section, we give an extension of Dubreil's Theorem
on the number of generators of ideals in a polynomial ring.
Specifically,  our generalization works for an ideal $I$ defining
a locally Cohen--Macaulay, equidimensional
subscheme $V$ of any codimension in $\Pn$, and
relates the number of generators of the defining ideal to
the length of certain Koszul homologies of
the cohomology of $V$.
The results in this section depend
crucially on the identification done in Section~1 of
the intersection modulo the product.

Finally, in Section~3, we give an extension of a surprising result
of Amasaki, \cite{Amasaki:structure},
showing a lower bound for the least degree of
a minimal generator of the ideal of a Buchsbaum subscheme.
Originally, Amasaki gave a bound in the case of Buchsbaum
curves in $\Bbb P^3$.  Easier proofs were subsequently given
by Geramita and Migliore in \cite{GM:generators},
based on a determination of the free resolution of the ideal
from a resolution of its deficiency module.  For Buchsbaum
codimension $2$ subschemes of $\Pn$ whose intermediate
cohomology vanishes, we are able to extend these considerations.

\section{When does $I \cap J = IJ$?}

Let $S= k[x_0, \dots, x_n]$ be a polynomial ring over the
algebraically closed field $k$.  Let $I$ and $J$ be ideals
defining subschemes $V$ and $Y$, respectively, of the
projective space $\Bbb P^n_k = \Pn$ over $k$.  In particular,
both $I$ and $J$ are homogeneous, saturated ideals.
In this section, we will derive a relationship between
the quotient module $(I \cap J)/IJ$ and the cohomology
of $V$, when $V$ and $Y$ meet in the expected dimension.

In general, if $V$ is a subscheme of $\Pn$, with saturated
homogeneous defining ideal $I = I_V$, the cohomology
modules of $V$ (or, less precisely, of $I$) are defined,
for $i=0, \dots, n-1$, by
$$
H^i_*({\cal I}_V) = H^i_*(V) = \bigoplus_j H^i(\Pn, {\cal I}_V(j)),
$$
where ${\cal I}_V = \widetilde{I_V}$ is the ideal sheaf of $V$.
These are all graded $S$-modules.  Moreover,
$H^0_*(V) = I_V$ and $H^i_*(V) = 0$ for $i > \dim V + 1$.
Usually, when $i = 1, \dots, \dim V$,
we will call $H^i_*(V)$ a deficiency module. This name
comes from the fact that the $H^i_*(V)$, $i = 1, \dots, \dim V$,
measure
the failure of $V$ to be an arithmetically Cohen--Macaulay
subscheme, since they vanish whenever $V$ is aCM.

We will also have need to use the cohomology of modules.
If $M$ is a (finitely generated) $S$-module,
let $\widetilde{M}$ be its sheafification.  Then, exactly as in
the case of ideal sheaves, we define the cohomology module of $\widetilde{M}$
to be
$$
H^i_*(\widetilde{M}) = \bigoplus_j H^i(\Pn, \widetilde{M}(j)).
$$
These are again graded $S$-modules.
We note here for future reference
that the cohomology modules of $M$ are related to the local
cohomology modules $H^i_\m(M)$ of $M$ with respect to the homogeneous
maximal ideal $\m$ as follows:
\begin{equation}\label{local-coh}
0 \rightarrow H^0_\m(M) \rightarrow M \rightarrow H^0_*(\widetilde{M})
	\rightarrow H^1_\m(M) \rightarrow 0
\end{equation}
$$
H^i_\m(M) \cong H^i_*(\widetilde{M})\mbox{\quad\quad for $i>1$.}
$$
See \cite[Chapter 0]{SV:buchsbaum} for a good discussion of
graded and local cohomology.

In this paper, we will sometimes require subschemes of $\Pn$ to
be locally Cohen--Macaulay and equidimensional.  This
is equivalent to saying that all the cohomology modules
have finite length, except of course for the top cohomology
$H^{d+1}_*(V)$, $d = \dim V$.
By Serre's vanishing theorem, this is again equivalent to
having $[H^i_*(V)]_j = 0$ for $j \ll 0$, and $1 \le i \le d$, since
in any case the cohomology modules vanish in high degrees.

Now, let $I$ and $J$ be as above, let $s = \projdim J$ be the
projective dimension of $J$, and write a minimal graded free resolution
of $J$ as follows:
\begin{equation}
\begin{array}{ccccccccccccccccc}
0 & \rightarrow & F_s &
 	\buildrel {\phi_s} \over \longrightarrow
	& \dots
	& {\buildrel {\phi_2} \over \rightarrow}
	& F_1&
	\buildrel {\phi_1} \over \longrightarrow & F_0 &
	\rightarrow & J & \rightarrow & 0.
\end{array}
\end{equation}
where $F_j = \bigoplus_i S(-a_{ji})$ are free modules.
For each $j = 1, \dots, s$, let $K_j$ be the $j$-th syzygy module, so that
there are short exact sequences
$$
0 \rightarrow K_{j+1} \buildrel {\psi_{j+1}} \over \longrightarrow
	F_j \buildrel {\eta_j} \over \longrightarrow K_j \rightarrow 0,
$$
where the maps $\psi_j$ and $\eta_j$ are the canonical inclusions and
projections, respectively. Note that for $j = s$, we have $K_s = F_s$,
$\psi_s = \phi_s$ and $\eta_s = id$.

For $S$-modules $M$ and $N$, and a map $f : M \to N$, we denote
by $f^i : H^i_*(\widetilde{M \otimes I}) \to H^i_*(\widetilde{N \otimes I})$
the map
induced on cohomology by $f \otimes id : M \otimes I \to N \otimes I$.

Our main technical result for this paper is the following Theorem.

\begin{thm}\label{main:technical}
Suppose the ideals $I$ and $J$ as above define
disjoint subschemes $V$ and $Y$, respectively.
Then for each $i \ge 1$, there are isomorphisms
$$
\ker \psi_i^1 \cong \Tor_i^S(S/I, S/J)
$$
and, for each $i,j \ge 1$, a long exact sequence
\begin{equation}\label{main:sequence}
0 \rightarrow \im \psi_{i+1}^j \rightarrow \ker \phi_i^j \rightarrow
	\ker \psi_i^j \rightarrow \ker \psi_{i+1}^{j+1} \rightarrow
	\coker \phi_i^j \rightarrow \coker \psi_i^j \rightarrow 0.
\nonumber
\end{equation}
\end{thm}

\begin{proof}  We remark that for $i > s$, both statements
are trivial, since then $\phi_i = \psi_i$ is the zero map.  Also,
if $j > \dim S/I$, then because $H^j_*(V) = 0$
again the second statement is trivial.

Now, let $\mu : I \otimes J \to IJ$ be
the natural surjection, and note that $\ker \mu \cong \Tor_2^S(S/I, S/J)$.
This follows, for instance, by tensoring
$$
0 \rightarrow I \rightarrow S \rightarrow S/I \rightarrow 0
$$
with $J$, comparing the resulting sequence with
$$
0 \rightarrow IJ \rightarrow J \rightarrow J/IJ \rightarrow 0
$$
via the multiplication map, and using that
$\Tor_1^S(S/I, J) \cong \Tor_2^S(S/I, S/J)$.  Note especially that
$\ker \mu$ has finite length since it is annihilated by $I+J$.
In particular, by sheafifying and taking cohomology of the
short exact sequence
$$
0 \rightarrow \ker \mu \rightarrow I \otimes J \rightarrow IJ
	\rightarrow 0,
$$
we see that $H^i_*(\widetilde{I \otimes J}) \cong H^i_*(\widetilde{IJ})$
for all $i \ge 0$.

Next, using the functorial map $M \to H^0_*(\widetilde{M})$ for any
$S$-module $M$, we get a commutative diagram
\begin{equation}\label{eq:1}
\begin{array}{ccccccccc}
&&&& 0 && 0 \\
&&&& \downarrow && \downarrow \\
&&&& H^0_\m(I \otimes J) & \rightarrow & H^0_\m(IJ) \\
&&&& \downarrow && \downarrow \\
0 &\rightarrow & \Tor_2^S(S/I, S/J) & \rightarrow & I \otimes J &
	{\buildrel {\mu} \over {\longrightarrow}} & IJ & \rightarrow & 0 \\
&&&& \downarrow && \downarrow \\
&& 0 & \rightarrow & H^0_*(\widetilde{I \otimes J}) &\rightarrow &
	H^0_*({\widetilde{IJ}}) & \rightarrow & 0 \\
&&&& \downarrow && \downarrow \\
&&&& H^1_\m(I \otimes J) & \rightarrow &H^1_\m(IJ) \\
&&&& \downarrow && \downarrow \\
&&&& 0 && 0
\end{array}
\end{equation}
But $IJ$ is an ideal, so $H^0_\m(IJ) = 0$.  Hence the kernel
of the map $I \otimes J \to H^0_*(\widetilde{I \otimes J})$
is $\Tor_2^S(S/I, S/J)$.

Now, with these preliminaries out of the way, we prove the
isomorphisms by induction on $i$.  For $i = 1$, tensor the exact
sequence
$$
0 \rightarrow K_1 {\buildrel {\psi_1} \over \longrightarrow }
	F_0 \rightarrow J \rightarrow 0
$$
by $I$.  This yields an exact sequence
$$
0 \rightarrow \Tor_1^S(I, J) \rightarrow K_1 \otimes I
	{\buildrel {\psi_1 \otimes 1} \over {\kindalong}}
	F_0 \otimes I \rightarrow J \otimes I \rightarrow 0.
$$
Since $\Tor_1^S(I, J) \cong \Tor_3^S(S/I, S/J)$ is annihilated
by $I+J$, in particular it has finite length.  Thus, taking cohomology
and comparing with the original sequence gives a diagram
\begin{equation}\nonumber
\setlength{\arraycolsep}{1pt}
\begin{array}{ccccccccccccc}
0 & \rightarrow & \Tor_3^S(S/I, S/J) & \rightarrow & K_1 \otimes I &
	\rightarrow & F_0 \otimes I & \rightarrow & J \otimes I &
	\rightarrow & 0 \\
&&&& \downarrow &&\downarrow&&\downarrow \\
&& 0 & \rightarrow & H^0_*(\widetilde{K_1 \otimes I}) & \rightarrow &
	H^0_*({\widetilde{F_0 \otimes I}}) & \rightarrow &
	H^0_*(\widetilde{J \otimes I}) & \rightarrow & \ker \psi_1^1 &
	\rightarrow & 0.
\end{array}
\end{equation}
Here, the middle vertical map is an isomorphism, since $I$ is
saturated.  Thus the snake lemma shows that there are exact sequences
$$
0 \rightarrow \Tor_3^S(S/I, S/J) \rightarrow K_1 \otimes I \rightarrow
	H^0_*(\widetilde{K_1 \otimes I}) \rightarrow \Tor_2^S(S/I, S/J)
	\rightarrow 0,
$$
and
$$
0 \rightarrow \Tor_2^S(S/I, S/J) \rightarrow J \otimes I \rightarrow
	H^0_*(\widetilde{J \otimes I}) \rightarrow \ker \psi_1^1
	\rightarrow 0.
$$
But from the above discussion, the last sequence implies the short exact
sequence
$$
0 \rightarrow IJ \rightarrow H^0_*(\widetilde{IJ}) \rightarrow \ker \psi_1^1
	\rightarrow 0.
$$
Since $I$ and $J$ define disjoint varieties, we have
$H^0_*(\widetilde{IJ}) = I \cap J$.  Thus, the above sequence
shows that
$$
\ker \psi_1^1 \cong {{I \cap J} \over {IJ}} \cong \Tor_1^S(S/I, S/J).
$$

By induction, we may assume that $\ker \psi_i^1 \cong \Tor_i^S(S/I, S/J)$,
and that there is an exact sequence
$$
0 \rightarrow \Tor_{i+2}^S(S/I, S/J) \rightarrow K_i \otimes I
	\rightarrow H^0_*(\widetilde{K_i \otimes I}) \rightarrow
	\Tor_{i+1}^S(S/I, S/J) \rightarrow 0.
$$
Tensoring the exact sequence
$$
0 \rightarrow K_{i+1} {\buildrel {\psi_{i+1}} \over \kindalong}
	F_i \rightarrow K_i \rightarrow 0
$$
with $I$ yields
$$
0 \rightarrow \Tor_1^S(K_i, I) \rightarrow K_{i+1} \otimes I
	\rightarrow F_i \otimes I \rightarrow K_i \otimes I \rightarrow 0.
$$
Here, $\Tor_1^S(K_i, I) \cong \Tor_{i+1}^S(I, J) \cong \Tor_{i+3}^S(S/I, S/J)$.
In particular, it is finite length.  Hence, taking cohomology and comparing
yields a diagram
\begin{equation}\nonumber
\setlength{\arraycolsep}{1pt}
\begin{array}{ccccccccccccc}
0 & \rightarrow & \Tor_{i+3}^S(S/I, S/J) & \rightarrow & K_{i+1} \otimes I &
	\rightarrow & F_i \otimes I & \rightarrow & K_i \otimes I &
	\rightarrow & 0 \\
&&&& \downarrow &&\downarrow&&\downarrow \\
&& 0 & \rightarrow & H^0_*(\widetilde{K_{i+1} \otimes I}) & \rightarrow &
	H^0_*({\widetilde{F_i \otimes I}}) & \rightarrow &
	H^0_*(\widetilde{K_i \otimes I}) & \rightarrow & \ker \psi_{i+1}^1 &
	\rightarrow & 0.
\end{array}
\end{equation}
But by the inductive hypothesis, we know the kernel and cokernel
of the right-hand vertical map.  Thus the snake lemma implies that
$$
\ker \psi_{i+1}^1 \cong \Tor_{i+1}^S(S/I, S/J),
$$
and that there is a long exact sequence
$$
0 \rightarrow \Tor_{i+3}^S(S/I, S/J) \rightarrow K_{i+1} \otimes I
	\rightarrow H^0_*(\widetilde{K_{i+1} \otimes I}) \rightarrow
	\Tor_{i+2}^S(S/I, S/J) \rightarrow 0,
$$
which finishes  the proof of the isomorphisms.

Next, we show that the long exact sequence exists.
Fix an $i \ge 1$.   Thus there is an exact sequence
$$
0 \rightarrow K_{i+1}
	{\buildrel {\psi_{i+1}} \over \kindalong} F_i
	{\buildrel {\eta_{i}} \over \kindalong} K_i \rightarrow 0.
$$
Tensor this sequence with $I$,  to obtain
$$
0 \rightarrow \Tor_1^S(K_i, I) \rightarrow K_{i+1} \otimes I
	{\buildrel {\psi_{i+1} \otimes 1} \over \kindalong} F_i \otimes I
	{\buildrel {\eta_{i} \otimes 1} \over \kindalong} K_i \otimes I
	\rightarrow 0,
$$
and note that
$\Tor_1^S(K_i , I) = \Tor_{i+3}^S(S/I, S/J)$, has finite length.
Thus, after sheafifying and taking cohomology, at the $j$-th stage
this yields isomorphisms
\begin{eqnarray*}
\ker \eta_i^j &\cong& \im \psi_{i+1}^j \\
\coker \eta_i^j &\cong& \ker \psi_{i+1}^{j+1}.
\end{eqnarray*}

Now, using the functoriality of tensor products and of cohomology,
we obtain a commutative square
\begin{equation}
\begin{array}{ccc}
H^j_*(\widetilde{F_i \otimes I})  & \widetilde{\kindalong} &
			H^j_*(\widetilde{F_i \otimes I}) \\
\mapdown{\eta_i^j} && \mapdown{\phi_i^j} \\
H^j_*(\widetilde{K_i \otimes I}) &
	{\buildrel  {\psi_i^j} \over {\longrightarrow}} &
		H^j_*(\widetilde{F_{i-1} \otimes I}).
\end{array}
\end{equation}
Applying the snake lemma to the columns, and  using the two isomorphisms
above shows that there is a sequence
$$
0 \rightarrow \im \psi_{i+1}^j \rightarrow \ker \phi_i^j \rightarrow
	\ker \psi_i^j \rightarrow \ker \psi_{i+1}^{j+1} \rightarrow
	\coker \phi_i^j \rightarrow \coker \psi_i^j \rightarrow 0,
$$
which is what we claimed.
\end{proof}

We note that this greatly extends the arguments in
\cite[Section 1]{Mig:submodules}.  The situation there was much
simpler in that it only considered the case that $J$ was
codimension $2$ and arithmetically Cohen--Macaulay (so most of
the terms in the sequence~(\ref{main:sequence}) vanish), and
only the case $i=j=1$ was studied, so it focused on
$\ker \psi_1^1 = (I \cap J)/IJ$.  Our
extension makes no assumptions on the Cohen--Macaulayness of $J$,
nor on its codimension.  Of course, our conclusion is much
more complicated, reflecting the fact that so much information
is encoded in the free resolution of $J$.

As an application of this technical result, in the next theorem
we give a proof of a statement due to Serre on when there
is an equality $I \cap J = IJ$.

\begin{thm}{\rm \cite[Corollaire, p. 143]{Serre}}
Suppose the ideals $I$ and $J$
define disjoint subschemes of ${\Bbb P}^n$.  Then $IJ = I \cap J$ if
and only if $\dim S/I + \dim S/J = \dim S$ and
both $S/I$ and $S/J$ are Cohen--Macaulay.
\end{thm}

\begin{proof}  Suppose first that $S/I$ and $S/J$ are
Cohen--Macaulay with $\dim S/I + \dim S/J = \dim S$.
Then by the Auslander--Buchsbaum formula,
$s = \projdim J = \dim S/I - 1$, and moreover $H^i_*(\widetilde{I}) = 0$
for $i = 1, \dots, s$.  In particular, $\ker \phi_i^i = 0 = \coker \phi_i^i$
for $i = 1, \dots, s$.  Since $\psi_s = \phi_s$, by reverse induction
the sequence (\ref{main:sequence}) with $j=i$ shows
that $\ker \psi_i^i = 0$ for $i = 1, \dots, s$.
Thus $(I \cap J)/IJ = \ker \psi_1^1 = 0$.



Conversely, since the subschemes defined by $I$ and $J$ are disjoint,
we have $\dim S/I + \dim S/J \le \dim S$  Hence
\begin{equation}\label{ineqs}
\dim S/I \le \dim S - \dim S/J \leq \dim S - \depth S/J = s+1
\end{equation}
where the last equality is by the Auslander--Buchsbaum formula.  Now,
if $IJ = I \cap J$, then $\Tor_1^S(S/I, S/J) = 0$, and
so by rigidity, $\Tor_i^S(S/I, S/J) = 0$ for $i \ge 1$.
Hence the isomorphisms of Theorem~\ref{main:technical}
show that $\ker \phi_s^1 = 0$.  But this implies that
$H^1_*(I) = 0$.  Thus also $\ker \phi_i^1 = 0 = \coker \phi_i^1$
for all $i = 1, \dots, s$, and since $\ker \psi_1^1 = (I \cap J)/IJ = 0$,
the exact sequence (\ref{main:sequence}) with $j = 1$, $i = 1, \dots, s$
implies that $\ker \psi^2_i = 0$, for $i = 2, \dots, s$.
Since $\psi_s = \phi_s$, this shows $\ker \psi_s^2 = 0$ and
hence also $H^2_*(\widetilde{I}) = 0$.  Continuing inductively,
we see that $\ker \psi_i^j = 0$ for all $i$ and $j$ with $i \ge j$.
In particular, since $\psi_s = \phi_s$ we get that
$H^j_*(\widetilde{I}) = 0$ for $j = 1, \dots, s$.

Let $d = \dim S/I$.  We have seen that $d-1 \le s$.  If this
inequality were strict, then in particular $H^d_*(\widetilde{I}) = 0$,
which is impossible.  Hence we have $d-1 = s$ and $H^j_*(\widetilde{I}) = 0$
for $j = 1, \dots, d-1$;  that is $S/I$ is Cohen--Macaulay.
But furthermore, each of the inequalities in (\ref{ineqs}) is
actually in equality.  This shows both that $\dim S/J = \depth S/J$,
i.e., $S/J$ is Cohen--Macaulay, and that $\dim S/I + \dim S/J = \dim S$,
which finishes the proof.
\end{proof}

\section{An Extension of Dubreil's Theorem}

In this section, we wish to use the results of Section 1 to
extend a theorem of Dubreil on the number of generators of
certain ideals.  Let $\nu(I)$ denote the minimal number of
generators of $I$, and let $\alpha(I)$ denote the
least degree of a minimal generator.  In its most
basic form, Dubreil's Theorem states:

\begin{thm}  Let $I$ be a homogeneous ideal of $k[x,y]$.  Then
$\nu(I) \le \alpha(I) + 1$.
\end{thm}

See \cite{DGM:Dubreil} for a proof of this;  note however that it
is essentially a consequence of the Hilbert--Burch theorem.
Dubreil's theorem is easily extended to the case that $I$
is a codimension $2$ arithmetically Cohen--Macaulay ideal in
any polynomial ring $k[x_0, \dots, x_n]$; again, see \cite{DGM:Dubreil}
for the details.  On the other hand, when $I$ is not arithmetically
Cohen--Macaulay, or when $I$ is not codimension $2$, not much is
known in this direction.  However, in the case of an ideal defining a
subscheme of $\Bbb P^3$, the following theorem of Migliore
shows that the general case will involve the cohomology of the
subscheme.

\begin{thm}{\rm \cite[Corollary 3.3]{Mig:submodules}}
Suppose $I$ defines a subscheme $V$ of $\Bbb P^3$,
of codimension at least $2$.  Let $A = (L_1, L_2)$ be the
complete intersection of two general linear forms, and
let $K_A$ denote the submodule of $H^1_*(V)$
annihilated by $A$.  Then
$$
\nu(I) \le \alpha(I) + 1 + \nu(K_A).
$$
\end{thm}

We note in particular that this formula is valid both for
the case that $V$ is codimension $2$, not necessarily
arithmetically Cohen--Macaulay, and the case that $V$ is
codimension $3$.  In the latter case, even though
$H^1_*(V)$ is not finitely generated, we still have
that at least $K_A$ is finitely generated
(see \cite[Theorem 2.1]{Mig:submodules} or our Lemma~\ref{finite-length}),
so the theorem still has useful content.  Furthermore, in
case $I$ defines an arithmetically Buchsbaum curve, so
that $H^1_*(V)$ is a $k$-vector space, $K_A = H^1_*(V)$
and $\nu(K_A) = \dim_k H^1_*(V)$.  Thus the Buchsbaum
case is particularly easy to calculate in examples.  In this
special case, the bound can be obtained from \cite{Amasaki:structure}.

In this section, we will give a generalization of
Dubreil's Theorem to ideals defining subschemes of $\Pn$ of
arbitrary codimension.  As an easy consequence, we recover by
our methods the above two theorems, and also part of a result of
Chang, \cite{Chang:charac}, on
the number of generators of an ideal defining a Buchsbaum
codimension $2$ subscheme of $\Pn$, which again seems
to be the best understood case.  Our generalization is a
corollary to the technical statement Theorem~\ref{main:technical}
in Section~1, underscoring the usefulness of identifying
the difference between intersections and products.

Our generalization will be based on the Koszul homologies
of the cohomology modules of an ideal $I$ defining a subscheme
of projective space.  As such, we will make some general
remarks concerning Koszul homology.  These comments are
basic, and can be found, for instance, in \cite{Mat}.
We first set the notation.  If $R$ is a ring, and
$y_1, \dots, y_s$ elements of $R$, we let ${\Bbb K}((y_1, \dots, y_s);R)$
denote the Koszul complex with respect to $y_1, \dots, y_s$.
If $M$ is an $R$-module, put
${\Bbb K}((y_1, \dots, y_s); M) = {\Bbb K}((y_1, \dots, y_s);R) \otimes M$,
the Koszul complex on $M$ with respect to $y_1, \dots, y_s$.
Set ${\Bbb H}_i((y_1, \dots, y_s); M)$ to be the $i$-th homology
module of ${\Bbb K}((y_1, \dots, y_s); M)$;  this is the
Koszul homology on $M$ with respect to $y_1, \dots, y_s$.
We will need the following facts:

\begin{remark}\label{Koszul}

\begin{enumerate}
\item  If $y_1, \dots, y_n$ forms a regular sequence on $M$,
then ${\Bbb K}((y_1, \dots, y_n); M)$ is acyclic.
\item  Let $J = (y_1, \dots, y_s)$.  Then for each $i = 0, \dots, s$,
$J \subseteq \ann {\Bbb H}_i((y_1, \dots, y_s);M)$.
\item  Suppose there is a short exact sequence of $R$-modules
$$
0 \rightarrow M_1 \rightarrow M_2 \rightarrow M_3 \rightarrow 0.
$$
Then there is a long exact sequence on Koszul homology
$$
\begin{array}{c}
\cdots \rightarrow
	{\Bbb H}_{i+1}((y_1, \dots, y_n); M_3) \rightarrow
	{\Bbb H}_i((y_1, \dots, y_n); M_1)
	\hspace{1.5in} \\
	\hspace{1.5in} \rightarrow
	{\Bbb H}_i((y_1, \dots, y_n); M_2) \rightarrow
	{\Bbb H}_i((y_1, \dots, y_n); M_3) \rightarrow \cdots.
\end{array}
$$
\item  For each $i = 0, \dots, s$, there is an isomorphism
$$
{\Bbb H}_i((0, y_2, \dots, y_s); M) \cong {\Bbb H}_i((y_2, \dots, y_s); M)
				\oplus {\Bbb H}_{i-1}((y_2, \dots, y_s); M).
$$
\end{enumerate}
\end{remark}

Throughout this section, let $I$ be the saturated defining
ideal of a locally Cohen--Macaulay, equidimensional subscheme
$V$ of ${\Bbb P}^n$; put $d=\dim V$. Let $J = (L_1, \dots, L_{n-1})$ be
the complete intersection of $n-1$ general linear forms.  In
particular, ${\Bbb K}((L_1, \dots, L_{n-1}); R)$ is a
free resolution of $S/J$.

Recall that the highest non-zero cohomology module
$H^{d+1}_*(\widetilde{I})$ is never finitely generated,
when $d \ge 0$.  In the next result, we show that nonetheless,
most of its Koszul homologies are finitely generated.
As general notation, for a module $M$ and an ideal $A$,
let $M_A$ denote the submodule of $M$ which is annihilated
by $A$; that is, $M_A = (0 :_M A)$.

\begin{prop}\label{finite-length}  Let $I$ and $J = (L_1, \dots, L_{n-1})$
be as above. Then
the Koszul homology
${\Bbb H}_i((L_1, \dots, L_{n-1});H^{d+1}_*(\widetilde{I}))$
is finitely generated for each $i \ge d+2$.  In particular,
if $d \ge 0$, the Koszul homology has finite length.
\end{prop}

\begin{proof} By changing coordinates if necessary, we may,
without loss of generality,
assume that $L_i = x_{i-1}.$  We will prove the proposition
by using induction on $d$.  If $d = -1$, i.e., $I$
defines the empty subscheme, then
$H^0_*(\widetilde{I}) = I$ is already finitely generated,
so each of its Koszul homologies is also finitely generated.

Suppose $d \ge 0$.  The exact sequence
$$
0 \rightarrow I {\buildrel {x_0} \over \longrightarrow} I
	\rightarrow I/x_0 I \rightarrow 0
$$
induces the long exact sequence on cohomology
\begin{equation}\label{endcohomology}
\begin{array}{rcl}
0 \rightarrow A \rightarrow H^{d}_*(\widetilde{I/x_0 I})
 & \verylong &
H^{d+1}_*(\tilde{I}) {\buildrel {x_0} \over \longrightarrow}
H^{d+1}_*(\tilde{I}) \rightarrow 0 \\
& \searrow \hfill \nearrow \\
& H^{d+1}_*(\tilde{I})_{(x_0)} \\
& \nearrow \hfill \searrow \\
0 && 0
\end{array}
\end{equation}
where $A = H^{d}_*(\widetilde{I}) / x_0 H^{d}_*(\widetilde{I})$.
Note in particular that $A$ is finitely generated, since it
is the quotient of two finitely generated modules.

{}From the right-hand part of this sequence, we obtain a long exact
sequence of Koszul homology (see Remark~\ref{Koszul})
$$
\begin{array}{c}
\cdots {\buildrel {x_0} \over \longrightarrow}
	{\Bbb H}_{i+1}((x_0, \dots, x_{n-2}); H^{d+1}_*(\widetilde{I}))
	\rightarrow
	{\Bbb H}_i((x_0, \dots, x_{n-2}); H^{d+1}_*(\widetilde{I})_{(x_0)})
	\hspace*{1in}\\
	\hspace*{1in}\rightarrow
	{\Bbb H}_i((x_0, \dots, x_{n-2}); H^{d+1}_*(\widetilde{I}))
	{\buildrel {x_0} \over \longrightarrow}
	{\Bbb H}_i((x_0, \dots, x_{n-2}); H^{d+1}_*(\widetilde{I}))
	\rightarrow \cdots.
\end{array}
$$
But for any module $M$, ${\Bbb H}_i((x_0, \dots, x_{n-2}); M)$
is annihilated by $x_0$, so the long exact sequence breaks into
short exact sequences
\begin{center}
\makebox[\textwidth][l]{$\qquad 0 \longrightarrow
	{\Bbb H}_{i+1}((x_0, \dots, x_{n-2}); H^{d+1}_*(\widetilde{I}))
	\longrightarrow
	{\Bbb H}_i((x_0, \dots, x_{n-2}); H^{d+1}_*(\widetilde{I})_{(x_0)})
	\longrightarrow$}
\makebox[\textwidth][r]{${\Bbb H}_i((x_0, \dots, x_{n-2});
H^{d+1}_*(\widetilde{I}))
	\longrightarrow 0.\qquad$}
\end{center}
Thus to show that ${\Bbb H}_i((x_0, \dots, x_{n-2}); H^{d+1}_*(\widetilde{I}))$
is finitely generated for $i \ge d+2$, it will suffice to show
that ${\Bbb H}_i((x_0, \dots, x_{n-2}); H^{d+1}_*(\widetilde{I})_{(x_0)})$
is finitely generated for $i \ge d+2$.  However, $x_0$ kills
$H^{d+1}_*(\widetilde{I})_{(x_0)}$, and so
{\setlength{\arraycolsep}{0pt}
\begin{eqnarray*}
{\Bbb H}_i((x_0, \dots, x_{n-2}); H^{d+1}_*(\widetilde{I})_{(x_0)})
	=&& {\Bbb H}_i((0, x_1, \dots, x_{n-2}); H^{d+1}_*(\widetilde{I})_{(x_0)}) \\
	= {\Bbb H}_i((x_1, \dots, &&x_{n-2}); H^{d+1}_*(\widetilde{I})_{(x_0)})
		\oplus {\Bbb H}_{i-1}((x_1, \dots, x_{n-2});
H^{d+1}_*(\widetilde{I})_{(x_0)}),
\end{eqnarray*}}
and we can calculate this over $R = S/(x_0) = k[x_1, \dots, x_n]$.

Now, the left-hand part of (\ref{endcohomology}) yields a long
exact sequence of Koszul homology
\begin{equation}\label{sequence1}
\begin{array}{c}
\cdots \rightarrow
	{\Bbb H}_{j+1}((x_1, \dots, x_{n-2}); H^{d+1}_*(\widetilde{I})_{(x_0)})
	\rightarrow {\Bbb H}_j((x_1, \dots, x_{n-2}); A)
	\rightarrow
	\hspace*{1in} \\
	\hspace*{1in}
	{\Bbb H}_j((x_1, \dots, x_{n-2}); H^d_*(\widetilde{I/x_0 I}))
	\rightarrow
	{\Bbb H}_j((x_1, \dots, x_{n-2}); H^{d+1}_*(\widetilde{I})_{(x_0)})
	\rightarrow\cdots.
\end{array}
\end{equation}
Here, the saturation of $I/x_0 I \subseteq R$ defines a subscheme
$\overline{V}$ of ${\Bbb P}^{n-1}$, with $\dim \overline{V} = d-1$.
Thus, by the induction hypothesis,
${\Bbb H}_j(H^d_*(\widetilde{I/ x_0 I}))$ is finitely
generated for each $j \ge d+1$. In particular, since $i \ge d+2$,
this is true for $j = i, i-1$.  Also, since $A$ is finitely
generated, all of its Koszul homologies are also finitely generated.
But then (\ref{sequence1})
shows that both
${\Bbb H}_i((x_1, \dots, x_{n-2}); H^{d+1}_*(\widetilde{I})_{(x_0)})$
and ${\Bbb H}_{i-1}((x_1, \dots, x_{n-2}); H^{d+1}_*(\widetilde{I})_{(x_0)})$
are finitely generated.  This implies that
${\Bbb H}_i((x_0, \dots, x_{n-2}); H^{d+1}_*(\widetilde{I}))$
is finitely generated, which finishes the proof of the first
statement.

For the second statement, recall that the Serre vanishing theorem
says that $H^{d+1}(\widetilde{I}(t))$ vanishes for large $t$, and
hence the Koszul homology ${\Bbb H}_i((L_1, \dots, L_{n-1});
H^{d+1}_*(\widetilde{I}))$
also vanishes in high degrees.  But since it is finitely generated,
it must also vanish in low degrees, and we can conclude that
it must have finite length.
\end{proof}

\begin{thm}\label{extended-Dubreil}
Suppose $I$ defines a locally Cohen--Macaulay, equidimensional
subscheme $V$ of dimension $d$ of $\Pn$.
Let $J = (L_1, \dots, L_{n-1})$ be generated by $n-1$ general linear
forms, and let
$\Bbb H_i((L_1, \dots, L_{n-1}); H^j_*(\widetilde{I}))$
be the Koszul homologies of
$H^j_*(\widetilde{I})$ with respect to $J$.  Then
\begin{equation}\label{dubreil:formula}
\nu(I) \le \alpha(I) + 1 + \sum_{i=1}^{n-2}
	   \dim \Bbb H_{i+1}((L_1, \dots, L_{n-1});H^i_*(V)).
\end{equation}
\end{thm}

\begin{proof}  Note that we have an exact sequence
$$
0 \rightarrow {{I \cap J} \over {IJ}} \rightarrow {I \over {IJ}}
	\rightarrow {{I+J} \over J} \rightarrow 0.
$$
Hence,
$\nu(I) = \nu(I/IJ) \le \nu({{I+J}\over J}) + \nu({{I \cap J} \over IJ})$.
Now, since $(I+J)/J$ is an ideal in $S/J$, which is a polynomial
ring in two variables, Dubreil's Theorem applies, and says
that $\nu((I+J)/J) \le \alpha((I+J)/J) + 1$.  Since $J$ is
generated by general linear forms, $\alpha((I+J)/J) = \alpha(I)$,
and so $\nu((I+J)/J) \le \alpha(I) + 1$.

Thus it only remains to estimate $\nu((I \cap J)/IJ)$.
Note that $I$ and $J$ define disjoint schemes, since
$J$ is generated by general linear forms, and that
the Koszul complex ${\Bbb K}((L_1, \dots, L_{n-1});S)$ is
a free resolution of $S/J$.  In particular,
the isomorphisms and exact sequences of Theorem~\ref{main:technical}
hold.
For each $i = 1, \dots, n-1$,
let $P_i = \ker \phi_i^i/\im \psi_{i+1}^i$.  Note then, that
$\Bbb H_{i+1}(H^i_*(V))$ naturally maps surjectively onto $P_i$,
since $\im \psi_{i+1}^i$ contains $\im \phi_{i+1}^i$.  Thus,
in particular, $\dim P_i \le \dim \Bbb H_{i+1}(H^i_*(V))$.  Hence
it follows from Theorem~\ref{main:technical} that
\begin{eqnarray}
\nu((I \cap J) / IJ) \le \dim ((I \cap J)/IJ) &\le& \dim \Bbb H_2(H^1_*(V))
							+ \dim \ker\psi_2^2
					\nonumber \\
	&\le& \dim \Bbb H_2(H^1_*(V)) + \dim \Bbb H_3(H^2_*(V)) +
				\dim \ker \psi_3^3 \nonumber \\
	&\vdots & \nonumber\\
	&\le& \sum_{i=1}^{n-2} \dim \Bbb H_{i+1}(H^i_*(V)). \nonumber
\end{eqnarray}

\end{proof}

\begin{remark}  Since $H^i_*(V) = 0$ whenever $i > \dim V + 1$,
many of the terms in the formula~(\ref{dubreil:formula}) vanish.
For instance, if $V$ is a curve in ${\Bbb P}^5$, there are only
two terms coming from the cohomology of $V$.
\end{remark}

\begin{remark}  In general, for a finite length graded module,
$\nu(M)$ is much less than $\dim M$, and so we would like to be able
to replace ``$\dim$'' by ``$\nu$'' throughout in the above formula.
However, counting minimal generators is much more difficult in
general than counting vector space dimensions.
\end{remark}

One important case in which we can replace ``$\dim$'' by ``$\nu$'' is
when all but the top cohomology of $V$ is annihilated by the
maximal ideal.  Recall the definition:

\begin{definition}  A subscheme $V$ of $\Pn$ of dimension $d$
is said to be {\em arithmetically Buchsbaum} if $H^i_*(V)$
is annihilated by the maximal ideal for each $i=1, \dots, d$,
and if for each general linear subspace $Y$ of $\Pn$,
the cohomology $H^i_*(V \cap Y)$ is annihilated by the maximal ideal
for $i = 1, \dots, \dim V \cap Y$.
\end{definition}

Note that the condition on linear subspaces of $V$ is required,
for there are examples of subschemes $V$ whose cohomology
is annihilated by the maximal ideal which have
hypersurface sections whose cohomology is not annihilated
by the maximal ideal;  see for instance, \cite{Miy:graded}.
In general, we will not require the
full strength of this definition, only that
the cohomologies are annihilated by the maximal ideal.
Such subschemes are called {\em quasi--Buchsbaum}.

The next corollary was obtained for the case $n=3$ in \cite{Mig:submodules},
and a somewhat better bound is stated in \cite{Chang:charac} for
arithmetically Buchsbaum subschemes of ${\Bbb P}^n$.

\begin{cor}  Let $I$ define a codimension $2$ subscheme
of $\Pn$ which is quasi-Buchsbaum.  Then
$$
\nu(I) \le \alpha(I) + 1 +
	\sum_{i=1}^{n-2} {{n-1} \choose {i+1}} \dim_k H^i_*(V).
$$
\end{cor}

\begin{proof}  We  only have to note that since $H^i_*(V)$ is
annihilated by the maximal ideal for $i = 1, \dots, n-1$, and since
every non-zero entry in a matrix representation for $\phi_i$
is a linear form, then for each $i = 1, \dots, n-1$,
$\Bbb H_{i+1}(H^i_*(V)) = \ker \phi_i^i$ is a
direct sum of ${{n-1} \choose {i+1}}$
copies of (twists of) $H^i_*(V)$.
\end{proof}

More generally, we can apply the same kind of analysis to quasi-Buchsbaum
subschemes of arbitrary codimension, except that we now have to
consider the top cohomology as well.

\begin{cor}\label{number}
Let $V$ be a $d$-dimensional quasi-Buchsbaum subscheme of $\Pn$,
defined by the saturated ideal $I$. Let $J = (L_1, \dots, L_{n-1})$
be generated by $n-1$ general linear forms, and let
$\Bbb H_i(H^j_*(V))$ be the Koszul homologies of $H^j_*(V)$
with respect to $J$.  Then
$$
\nu(I) \le \alpha(I) + 1 +
 	\sum_{i=1}^{d} {{n-1} \choose {i+1}} \dim_k H^i_*(V)
		+ \dim \Bbb H_{d+2}(H^{d+1}_*(V)).
$$
\end{cor}

\begin{proof}  Again, since $H^i_*(V)$ is annihilated by every
linear form for $i = 1, \dots, d$, then $\Bbb H_{i+1}(H^i_*(V))$
is just a direct sum of ${{n-1} \choose {i+1}}$ copies of
$H^i_*(V)$.
\end{proof}

\begin{remark}  By Lemma~\ref{finite-length}, even though
$H^{d+1}_*(V)$ does not have finite length, the Koszul
homology $\Bbb H_{d+2}(H^{d+1}_*(V))$ does have
finite length, and so this corollary really does
give a finite bound on the number of generators of $I$.
\end{remark}

\begin{example} Unfortunately, the bound in
Theorem~\ref{extended-Dubreil} does not seem to be very sharp.
For example, in ${\Bbb P}^4$, let $V$ be the union of a conic and a line not
meeting the plane of the conic.
Then $V$ is arithmetically Buchsbaum,
and $\dim_k H^1_*(V) = 1$.  Also, $\alpha(I_V) = 2$,
and $\nu(I_V) = 7$.  However, if $J$ is generated by three
general  linear forms, then $H^2_*(V)_J$ is at least $2$,
as can be seen, for instance, by a calculation using
the inductive procedure in Proposition~\ref{finite-length}.
So $\alpha(I_V) + 1 + 3\dim_k H^1_*(V) + \dim H^2_*(V)_J \geq 8$,
but $\nu(I_V) = 7$.
\end{example}

\begin{remark}  During the final preparation of this paper,
the authors received the preprint \cite{Hoa} of Hoa,
which contains bounds on the number of generators of
an ideal based in part on the cohomology of the ideal, but
involving different invariants of the ideal than
our bounds. Neither Hoa's bounds nor our bounds seem to be
particularly sharp in general.
\end{remark}

\section{On the Least Degree of Surfaces Containing Certain Buchsbaum
		Subschemes}

In this section, we want to use the bound given in Section~2 to
extend a result of Amasaki on
the minimal degree of the minimal generators of an ideal $I$
defining a codimension $2$ Buchsbaum subscheme of $\Pn$.
In \cite{Amasaki:structure}, Amasaki showed that if $C$ is
a Buchsbaum curve in $\Bbb P^3$, and if $N = \dim_k H^1_*(C)$
is the Buchsbaum invariant of $C$, then $\alpha(I) \ge 2N$.
A different proof was subsequently given in \cite{GM:generators}
based on combining the upper bound estimate for $\nu(I)$ of
Corollary~\ref{number} in the case of curves in $\Bbb P^3$
together with a lower bound estimate coming from a determination
of the free resolution of the ideal $I$ from a free resolution
of $H^1_*(C)$.  Also, Chang extended Amasaki's bound to
a codimension $2$ Buchsbaum subscheme of any $\Pn$ in \cite{Chang:charac},
based on a structure theorem for the locally free resolution of
the ideal sheaf associated to the subscheme.

Here, we would like to use our methods to give a different
proof of Amasaki's bound for certain codimension $2$ subschemes
of $\Pn$.  Specifically, we will give a lower bound for $\alpha(I)$
in terms of $H^1_*(V)$ for a codimension $2$ subscheme of $V$
for which $H^1_*(V)$ is annihilated by the maximal ideal,
and $H^i_*(V) = 0$ for $i = 2, \dots, \dim V$.
Note that these quasi-Buchsbaum schemes are in fact Buchsbaum, since
if $H$ is a general hyperplane defined by a linear form $L$,
the from the standard exact sequence
$$
0 \rightarrow I {\buildrel {\times L} \over {\longrightarrow}} I
	\rightarrow I/LI \rightarrow 0,
$$
it is easy to see that $H^i_*(V \cap H) \cong H^i_*(V)$
for $1 \le i \le \dim V - 1$.
Our method of proof will be to follow the lines of
\cite{GM:generators}.  That is, we will combine our
upper bound estimate of $\nu(I)$ with a lower bound
estimate for $\nu(I)$ based on the free resolution of $I$.

We begin with an extension of a result in \cite{Rao}.

\begin{prop}  Suppose $V$ is a codimension $2$ subscheme
of $\Pn$, such that $H^i_*(V) = 0$ for all $i = 2, \dots, \dim V$.
Let
$$
0 \rightarrow L_{n+1} {\buildrel {\sigma_{n+1}} \over \longrightarrow} L_n
	{\buildrel {\sigma_n} \over \longrightarrow} L_{n-1} \rightarrow \dots
	{\buildrel {\sigma_1} \over \longrightarrow} L_0
	\rightarrow H^1_*(V) \rightarrow 0
$$
be the minimal free resolution of the finite length module $H^1_*(V)$.
Then the saturated defining ideal $I = I_V$ of $V$ has a minimal
free resolution
$$
0 \rightarrow L_{n+1} {\buildrel {\sigma_{n+1}} \over \longrightarrow}
	L_n {\buildrel {\sigma_{n}} \over \longrightarrow}
	L_{n-1} {\buildrel_{\sigma_{n-1}} \over \longrightarrow}
	\dots {\buildrel {\sigma_4} \over \longrightarrow}
	L_4 {\buildrel {{[\sigma_3\; 0]}} \over \kindalong}
	L_3 \oplus \bigoplus_{i=1}^r S(-b_i)
	\rightarrow \oplus_{j=1}^p S(-a_j)
	\rightarrow I \rightarrow 0,
$$
for some $r \ge 0$, and $p = \nu(I)$.
\end{prop}

\begin{proof}  Write a minimal free resolution of $S/I$ as
follows:
$$
0 \rightarrow F_n {\buildrel {\phi_n} \over \longrightarrow} F_{n-1}
	{\buildrel {\phi_{n-1}} \over \longrightarrow} \dots
	{\buildrel {\phi_{2}} \over \longrightarrow} F_1
	\rightarrow S \rightarrow S/I \rightarrow 0,
$$
and let $E_i$ be the $i$-th syzygy module.  Sheafifying and dualizing
the short exact sequence
$$
0 \rightarrow F_n {\buildrel {\phi_{n-1}} \over \longrightarrow}
	F_{n-1} \rightarrow E_{n-1} \rightarrow 0
$$
yields the exact sequence
$$
0 \rightarrow {\cal E}_{n-1}^\vee \rightarrow {\cal F}_{n-1}^\vee
	{\buildrel {\phi_{n-1}^\vee} \over \kindalong}
	{\cal F}_n^\vee \rightarrow 0.
$$
Taking cohomology then gives a sequence
$$
0 \rightarrow H^0_*({\cal E}_{n-1}^\vee) \rightarrow
	F_{n-1}^\vee {\buildrel {\phi_{n-1}^\vee} \over \longrightarrow}
	F_n^\vee \rightarrow
	H^1_*({\cal E}_{n-1}^\vee) \rightarrow 0.
$$
Note, though, that
$H^1_*({\cal E}_{n-1}^\vee) \cong \Ext_S^{n+1}(H^1_*(V), S))$.

Next, for each $i = 2, \dots, n-2$, consider the sequence
$$
0 \rightarrow E_{i+1} \rightarrow F_{i} \rightarrow E_{i} \rightarrow 0.
$$
Sheafifying, dualizing and taking cohomology gives an exact sequence
$$
0 \rightarrow H^0_*({\cal E}_i^\vee) \rightarrow F_{i}^\vee
	\rightarrow H^0({\cal E}_{i+1}^\vee) \rightarrow
	H^1_*({\cal E}_i^\vee) \rightarrow 0.
$$
But $H^1_*({\cal E}_i^\vee) = \Ext_S^{n+1}(H^{n-i}_*(V), S) = 0$,
by assumption.  Hence, we can paste together all these exact sequences
to get a long exact sequence
$$
F_2^\vee {\buildrel {\phi_{3}^\vee} \over \kindalong}
	F_3^\vee {\buildrel {\phi_{4}^\vee} \over \kindalong}
	\dots {\buildrel {\phi_{n}^\vee} \over \kindalong}
	\Ext_S^{n+1}(H^1_*(V), S) \rightarrow 0.
$$
However, a minimal free resolution of $\Ext_S^{n+1}(H^1_*(V), S)$
is given by just dualizing the resolution of $H^1_*(V)$, and
so we see that $F_i = L_{i+1}$  and
$\phi_i = \sigma_{i+1}$ for $i = 3, \dots, n$, and
$F_2 = L_3 \oplus \bigoplus_{i=1}^r S$ and $\phi_2 = [\sigma_3\; 0]$,
for some $r \ge 0$.  This finishes the proof.
\end{proof}

\begin{cor}  Let $V$ be as in the previous proposition, and
let $I$ be its saturated defining ideal.  Then
$\nu(I) \ge 1 + \sum_{i=3}^{n+1} (-1)^i \rank L_i$.
\end{cor}

\begin{proof}  With the notation as in the statement of the Proposition,
we have
$$
\nu(I) = p = 1 + \rank L_3 + r -\rank L_4 + \rank L_5 + \dots
	\ge 1 + \sum_{i=3}^{n+1} (-1)^i \rank L_i.
$$
\end{proof}

\begin{cor}  In addition to the assumptions of the Proposition,
suppose that $H^1_*(V)$ is annihilated by the maximal ideal.
Let $N = \dim_k H^1_*(V)$. Then $\alpha(I) \ge (n-2)N$.
\end{cor}

\begin{proof}  A minimal free resolution of $H^1_*(V)$ is
just a direct sum of $N$ copies of the Koszul complex
resolving $k = S/\m$.  Thus, $\rank L_i = N{{n+1} \choose i}$.
By the previous Corollary and Corollary~\ref{number},
we have
$$
1 + \sum_{i=3}^{n+1} (-1)^i N{{n+1} \choose i} \le 1 + \alpha(I)
							+ (n+1)N.
$$
A simple arithmetic calculation then reduces this
to $(n-2)N \le \alpha(I)$, as claimed.
\end{proof}

\vskip.5truein

\makebox[\textwidth]{\hskip.5in \parbox{3in}{Department of Mathematics\\
					    Florida State University \\
					    Tallahassee \\
					    Florida, 32306 \\
					    U.S.A}
		     \hfill     \parbox{3in}{Department of Mathematics \\
					    University of Notre Dame \\
					    Notre Dame\\
				    	    Indiana, 46556 \\
					    U.S.A}
		     \hskip.5in}

\end{document}